\documentstyle[12pt]{article}
\newcommand{\text}{\mbox}
\newcommand{\hf}{\hspace{.4in}}
\begin {document}
\begin{center}
{\large\bf Some Remarks on the Best Approximation Rate\\ of
Certain Trigonometric Series}\\

        Songping Zhou \footnote{Supported in part by National and Zhejiang Provincial
Natural Science Foundations of China under grant numbers 1041001
and 101009 respectively.}and Ruijun Le
\end{center}

          \pagenumbering{arabic}
          \begin{quote}
          \small \bf ABSTRACT.
\rm  The main object of the present paper is to give a complete result regarding the best approximation rate of certain trigonometric series in general
 complex valued continuous function space under a new condition which gives an essential generalization to  $O$-regularly varying quasimonotonicity.   An application is present in Section 3.

          \end{quote}

\begin{center}
\end{center}

\begin{center}
\large\bf \S 1. Introduction
\end{center}

\vspace{3mm}

\rm As generalizations to monotonicity condition, various quasimonotonicity conditions including so-called $O$-regularly varying quasimonotonicity are introduced and applied to convergence of trigonometric (Fourier) series.      Interested readers could check references such as
Nurcombe [3], and Xie and Zhou [8] for uniform convergence and Stanojevic [4, 5], and Xie and Zhou [7] for $L^{1}$ convergence.

\rm Let $C_{2\pi}$ be the space of
 all real or complex valued continuous functions $f(x)$ of period $2\pi$ with  norm
$$\|f\| = \max_{-\infty<x<\infty}|f(x)|,$$
and $E_{n}(f)$ the best approximation by trigonometric polynomials of degree $n$. Given a trigonometric series $\sum_{k=-\infty}^{\infty}c_{k}e^{ikx}:=\lim\limits_{n\to\infty}\sum_{k=-n}^{n}c_{k}e^{ikx}$, write
$$f(x) = \sum_{k=-\infty}^{\infty}c_{k}e^{ikx}
=:\sum_{k=-\infty}^{\infty}\hat{f}(k)e^{ikx}$$
at any point $x$ where the series converges. Denote its $n$th partial sum $S_{n}(f,x)$ by
$$\sum_{k=-n}^{n}\hat{f}(k)e^{ikx}.$$

For a sequence $\{c_{n}\}_{n=0}^{\infty}$, let
$$\Delta c_{n} = c_{n}-c_{n+1}.$$

A non-decreasing positive sequence $\{R(n)\}_{n=1}^{\infty}$ is said to be $O$-regularly varying if\footnote{In some papers, this requirement is written as that for some $\lambda>1$, $\limsup\limits_{n\to\infty}R([\lambda n])/R(n)<\infty$. It is just a pure decoration.}
$$\limsup_{n\to\infty}\frac{R(2n)}{R(n)}<\infty.$$

A complex sequence $\{c_{n}\}_{n=0}^{\infty}$ is $O$-regularly varying quasimonotone in complex sense if for some $\theta_{0}\in [0,\pi/2)$ and some $O$-regularly varying sequence $\{R(n)\}$ the sequence
$$\Delta\frac{c_{n}}{R(n)}\in K(\theta_{0}):=\{z: |\mbox{\rm arg} z|\leq\theta_{0}\},\;\;n=1, 2,\cdots.$$

Evidently, if $\{c_{n}\}$ is a real sequence, then the $O$-regularly varying quasimonotonicity becomes
$$\Delta\frac{c_{n}}{R(n)}\geq 0,\;\;n=1, 2, \cdots,$$
which was used in many works to generalize the regularly varying quasimonotone condition and, in particular, the quasimonotone condition\footnote{A real sequence $\{b_{n}\}_{n=0}^{\infty}$
is defined to be quasimonotone if, for some $\alpha\geq 0$, the sequence $\{b_{n}/n^{\alpha}\}$ is non-increasing.}.

In 1992, without proof, Belov [1] announced the following result:

\vspace{3mm}

\bf Theorem B. \it Let $\{a_{n}\}$ and $\{b_{n}\}$ be quasimonotone sequences satisfying
$\sum_{n=0}^{\infty}a_{n}<\infty$ and $\lim\limits_{n\to\infty}nb_{n}=0$. Define
$$f(x)=\sum_{n=0}^{\infty}a_{n}\cos nx,\hspace{.2in} g(x)=\sum_{n=1}^{\infty}b_{n}\sin nx.$$
Then $f,g\in C_{2\pi}$ have the following estimates:
$$E_{n}(f)\approx \max_{1\leq k\leq n}ka_{n+k}+\sum_{k=2n+1}^{\infty}a_{k},$$
and
$$E_{n}(g)\approx \max_{k\geq 1}kb_{n+k}.$$

\vspace{3mm}

\rm It is essentially generalized recently by Xiao, Xie and Zhou [6] for general
 complex valued trigonometric series (thus a proof is given to Theorem B) as follows.

\vspace{3mm}

\bf Theorem XZ. \it
Let
$\{\hat{f}(n)\}_{n=0}^{\infty}$ and $\{\hat{f}(n)+\hat{f}(-n)\}_{n=1}^{\infty}$ be  complex $O$-regularly varying quasimonotone sequences, and
$$f(x)=\sum_{n=-\infty}^{\infty}\hat{f}(n)e^{inx}.$$
Then $f\in C_{2\pi}$ if and only if
$$\sum_{n=1}^{\infty}|\hat{f}(n)+\hat{f}(-n)|<\infty$$
and
$$\lim\limits_{n\to\infty}n\hat{f}(n)=0.$$
Furthermore, if $f\in C_{2\pi}$, then\footnote{As usual, $A_{n}\approx B_{n}$ means that there is a positive constant $C>0$ independent of $n$ such that $C^{-1}B_{n}\leq A_{n}\leq CB_{n}$.}
$$E_{n}(f)\approx
\max_{1\leq k\leq n}k\left(\left|\hat{f}(n+k)\right|
+\left|\hat{f}(-n-k)\right|\right)
+\max_{k\geq 2n+1}k\left|\hat{f}(k)
-\hat{f}(-k)\right|$$
$$+\sum_{k=2n+1}^{\infty}|\hat{f}(k)+\hat{f}(-k)|.$$

\rm We make a quick remark here. We notice that, in the complex valued function space, the assumption that both $\{\hat{f}(n)\}_{n=0}^{\infty}$ and $\{\hat{f}(-n)\}_{n=1}^{\infty}$ are  complex $O$-regularly varying quasimonotone is a convenient one, but is almost trivial, since it is almost the same as the condition in real case. Thus people usually use one side quasimonotonicity with some kind balance condition in considering those problems. Theorem XZ reflects this kind of thinking.

\vspace{3mm}

\rm Motivated by an idea in  Leindler [2], our work [11] introduced a new condition, and we present here a revised form as follows:

\vspace{3mm}

\bf Definition. \it Let \mbox{\bf c}$=\{c_{n}\}_{n=1}^{\infty}$ be a sequence satisfying $c_{n}\in K(\theta_{1})$ for some $\theta_{1}\in [0, 2\pi)$ and $n=1, 2, \cdots$. If there is a natural number $N_{0}$ such that
$$\sum_{n=m}^{2m}|\Delta c_{n}|\leq M(\mbox{\bf c})\max_{m\leq n<m+N_{0}}|c_{n}|\hspace{.4in}(1)$$
holds for all $m=1, 2, \cdots$, where  $M(\mbox{\bf c})$ indicates a positive constant only depending upon $\mbox{\bf c}$, then we say that the sequence $\mbox{\bf c}$ belongs to class \mbox{\rm GBV}.

\vspace{3mm}

\rm We recall the following results.

\vspace{3mm}

\bf Lemma 1. \it Suppose a complex sequence $\{c_{n}\}$ is $O$-regularly varying quasimonotone,
then there is a positive constant $M$ depending upon $\theta_{0}$ only such that
$$|c_{n}| \leq M\text{\rm Re}c_{n},\;\;n=1, 2,\cdots,$$
or in other words, $c_{n}\in K(\theta_{1})$ for some
$\theta_{1}\in [0, \pi/2)$ and $n=1, 2, \cdots$.

\vspace{3mm}

 \rm The argument exactly follows from Xie and Zhou [8, Lemma 1], and the condition $\lim\limits_{n\to\infty}c_{n}=0$ there can be cancelled however.
\vspace{3mm}

\bf Lemma 2. \it Let $\{c_{n}\}$ be any given complex $O$-regularly varying quasimonotone sequence. Then $\{c_{n}\}$ satisfies $(1)$ for $N_{0}=1$.

\vspace{3mm}

\rm The proof can be copied from Zhou and Le [11, Theorem 3] with omitting the condition $\lim\limits_{n\to\infty}c_{n}=0$ there.

\vspace{3mm}

We obviously see that, from Lemma 1 and Lemma 2, if $\mbox{\bf c}=\{c_{n}\}_{n=1}^{\infty}$ is any given complex $O$-regularly varying quasimonotone sequence, then $\{c_{n}\}\in \mbox{\rm GBV}$. On the other hand, since the converse is not true (see [11]), the class \mbox{\rm GBV} gives an essential and explicit generalization to the class of $O$-regularly varying quasimonotone sequences.

\rm From $\sum_{k=1}^{\infty}k^{-\alpha}\sin 2^{k}x=:\sum_{n=1}^{\infty}b_{n}\sin nx$, $\alpha>1$, we can clearly see that, for any $\epsilon>0$, $n^{\epsilon}b_{n}\to \infty$, $n\to\infty$, therefore, the condition (1), in general sense, cannot be further generalized.

Throughout the paper, $C$ denotes a positive constant (which is independent of $n$ and $x\in [0,2\pi]$) not necessarily the same at each occurrence. In some specific cases, we also use $M(\mbox{\bf c})$ to indicate a positive constant only depending upon the sequence $\mbox{\bf c}$.

\vspace{3mm}

\begin{center}
\large\bf \S 2. Results and Proofs
\end{center}

\vspace{3mm}

\rm We present the main result of the paper.

\vspace{3mm}

\bf Theorem 1. \it
Let
$\{\hat{f}(n)\}_{n=0}^{\infty}\in \mbox{\rm GBV}$ and $\{\hat{f}(n)+\hat{f}(-n)\}_{n=1}^{\infty}\in \mbox{\rm GBV}$, and
$$f(x)=\sum_{n=-\infty}^{\infty}\hat{f}(n)e^{inx}.$$
Then $f\in C_{2\pi}$ if and only if
$$\sum_{n=1}^{\infty}|\hat{f}(n)+\hat{f}(-n)|<\infty\hspace{.4in}(2)$$
and
$$\lim\limits_{n\to\infty}n\hat{f}(n)=0.\hspace{.4in}(3)$$
Furthermore, if $f\in C_{2\pi}$, then
$$E_{n}(f)\approx
\max_{1\leq k\leq n}k\left(\left|\hat{f}(n+k)\right|
+\left|\hat{f}(-n-k)\right|\right)
+\max_{k\geq 2n+1}k\left|\hat{f}(k)
-\hat{f}(-k)\right|$$
$$+\sum_{k=2n+1}^{\infty}|\hat{f}(k)+\hat{f}(-k)|.$$

\rm First we establish several lemmas.

\vspace{3mm}

\bf Lemma 3 (Xie and Zhou [8, Lemma 2]). \it Let $\{\hat{f}(n)\}$ satisfy
$$\hat{f}(n)+\hat{f}(-n) \in K(\theta_{1}),\;\;n=1,2,\cdots,$$
for some $\theta_{1}\in [0,\pi/2)$.
Then $f\in C_{2\pi}$ implies that
$$\sum_{n=1}^{\infty}|\hat{f}(n)+\hat{f}(-n)| < \infty.$$

\bf Lemma 4. \it Let
$\{\hat{f}(n)\}_{n=0}^{\infty}\in \mbox{\rm GBV}$ and $\{\hat{f}(n)+\hat{f}(-n)\}_{n=1}^{\infty}\in \mbox{\rm GBV}$.
Suppose $f\in C_{2\pi}$, then
$$\max_{k\geq 1}k|\hat{f}(\pm (n+k))|=O(E_{n}(f)).$$

\bf Proof. \rm Let $t_{n}^{*}(x)$ be the trigonometric polynomials of best approximation of degree $n$, then from an easy equality
$$\frac{1}{2\pi}\int_{-\pi}^{\pi}\left|e^{\pm i(n+1)x}\left(\sum_{k=0}^{N-1}e^{\pm ikx}\right)^{2}\right|dx=N$$
we get
$$\left|\frac{1}{2\pi}\int_{-\pi}^{\pi}(f(x)-t_{n}^{*}(x))e^{\pm i(n+1)x}\left(\sum_{k=0}^{N-1}e^{\pm ikx}\right)^{2}dx\right|\leq NE_{n}(f).$$
Also, the left integral within the absolute value symbols of the above inequality equals to
$\sum\limits_{k=1}^{N}\left(k\hat{f}(\mp (n+k))+(N-k)\hat{f}(\mp (n+N+k))\right)$, then
$$NE_{n}(f)\geq
\left|\sum\limits_{k=1}^{N}\left(k\hat{f}(\mp (n+k))+(N-k)\hat{f}(\mp (n+N+k))\right)\right|$$
$$\geq \mbox{\rm Re}\left(\sum\limits_{k=1}^{N}\left(k\hat{f}(\mp (n+k))+(N-k)\hat{f}(\mp (n+N+k))\right)\right).$$
Thus by the definition of \mbox{\rm GBV},
$$NE_{n}(f) \geq
\sum\limits_{k=1}^{N}\left(k\mbox{\rm Re}\hat{f}(n+k)+ (N-k) \mbox{\rm Re}\hat{f}(n+N+k)\right)$$
$$\hspace{2.6in}\geq
\sum\limits_{k=1}^{N}k\mbox{\rm Re}\hat{f}(n+k), \hf (4)$$
and
$$2NE_{n}(f)\geq
 \sum\limits_{k=1}^{N}k\left(\mbox{\rm Re}\hat{f}(n+k)+ \mbox{\rm Re}\hat{f}(-n-k)\right). \hf (5)$$
By the same argument as (4), we also have
$$2NE_{2n}(f)\geq
\sum\limits_{k=1}^{2N}k\mbox{\rm Re}\hat{f}(2n+k). \hf (4')$$

Fix sufficient large $n$, assume
$$\max_{N/2+jN_{0}\leq k<N/2+(j+1)N_{0}}|\hat{f}(n+k)|=| \hat{f}(n+k_{j})|,\;\;j=0, 1, \cdots, [N/(2N_{0})]-1,$$
$$\hspace{1.0in}n+N/2\leq n+N/2+jN_{0}\leq n+k_{j}< n+N/2+(j+1)N_{0}\leq n+N.$$
From condition $(1)$, we get for $0\leq j\leq [N/(2N_{0})]-1$,
$$|\hat{f}(n+N)|=\left|\sum_{k=n+N}^{2n+N+2jN_{0}-1}\Delta \hat{f}(k)+ \hat{f}(2n+N+2jN_{0})\right|$$
$$\leq \sum_{n+N/2+jN_{0}\leq k\leq 2n+N+2jN_{0}}|\Delta \hat{f}(k)|+|\hat{f}(2n+N+2jN_{0})|$$
$$\leq M(\mbox{\bf $\hat{f}$})(|\hat{f}(n+k_{j})|+|\hat{f}(2n+N+2jN_{0})|)$$
$$\leq M(\mbox{\bf $\hat{f}$}, \theta_{1})\left(\mbox{\rm Re}\hat{f}(n+k_{j})+\mbox{\rm Re}\hat{f}(2n+N+2jN_{0})\right),$$
therefore, together with (4) and $(4')$, we deduce that
$$3NE_{n}(f)\geq  NE_{n}(f) + 2NE_{2n}(f)\geq
\sum\limits_{k=1}^{N}k\mbox{\rm Re}\hat{f}(n+k)+ \sum\limits_{k=1}^{2N}k\mbox{\rm Re}\hat{f}(2n+k)$$
$$\geq\sum_{j=0}^{[N/(2N_{0})]-1} k_{j}\mbox{\rm Re} \hat{f}(n+k_{j})+\sum_{j=0}^{[N/(2N_{0})]-1}(N+2jN_{0})\mbox{\rm Re}\hat{f}(2n+N+2jN_{0})$$
$$\geq C\sum_{j=0}^{[N/(2N_{0})]-1} k_{j}\left(\mbox{\rm Re} \hat{f}(n+k_{j})+\mbox{\rm Re}\hat{f}(2n+N+2jN_{0})\right)$$
$$\geq C\sum_{j=0}^{[N/(2N_{0})]-1} k_{j}|\hat{f}(n+N)| \geq C|\hat{f}(n+N)|\sum_{j=0}^{[N/(2N_{0})]-1}(N/2+jN_{0})$$
$$\geq CN_{0}^{-1}N^{2}|\hat{f}(n+N)|.$$
Finally, we achieve that, for $N\geq 1$,
$$ N|\hat{f}(n+N)|\leq CE_{n}(f),$$
or in other words,
$$ \max_{k\geq 1}k|\hat{f}(n+k)|\leq CE_{n}(f).$$
At the same time, starting from (5),  since $\{\hat{f}(n)+\hat{f}(-n)\}_{n=1}^{\infty}\in \mbox{GBV}$, a similar argument leads to that
$$\max_{k\geq 1}k|\hat{f}(n+k)+\hat{f}(-n-k)|=O(E_{n}(f)).$$
Thus
$$\max_{k\geq 1}k|\hat{f}(-n-k)|\leq \max_{k\geq 1}k|\hat{f}(n+k)+\hat{f}(-n-k)|+\max_{k\geq 1}k|\hat{f}(n+k)|=O(E_{n}(f)).$$
Lemma 4 is completed.

As an important application, we write the following corollary of Lemma 4 as a theorem.

\vspace{3mm}

\bf Theorem 2.
\it
Let
$\{\hat{f}(n)\}_{n=0}^{\infty}\in \mbox{\rm GBV}$ and $\{\hat{f}(n)+\hat{f}(-n)\}_{n=1}^{\infty}\in \mbox{\rm GBV}$.
If $f\in C_{2\pi}$, then
$$\lim_{n\to\infty}n|\hat{f}(\pm n)|=0.$$

\vspace{3mm}

\bf Corollary.
\it
Let
$\{b_{n}\}_{n=1}^{\infty}\in\mbox{\rm GBV}$ be a real sequence\footnote{It means $b_{n}\geq 0$ and satisfies condition (1).}.
Suppose $g(x)=\sum_{n=1}^{\infty}b_{n}\sin nx\in C_{2\pi}$, then
$$\lim_{n\to\infty}nb_{n}=0.$$

\vspace{3mm}

\bf Lemma 5 (Xiao, Xie and Zhou [6, Lemma 4]). \it
Let $f\in C_{2\pi}$, $\hat{f}(n)+\hat{f}(-n)\in K(\theta_{1})$, $n=1,2,\cdots$, for some $0\leq\theta_{1}<\pi/2$. Then
$$\sum_{k=2n+1}^{\infty}|\hat{f}(k)+\hat{f}(-k)|=O(E_{n}(f)).$$

 \bf Lemma 6. \it
Let
$\{\hat{f}(n)\}_{n=0}^{\infty}\in \mbox{\rm GBV}$ and $\{\hat{f}(n)+\hat{f}(-n)\}_{n=1}^{\infty}\in \mbox{\rm GBV}$,
then
$$\left|\sum_{k=1}^{n}\hat{f}(\pm(n+k))\sin kx\right|=
O\left(\max_{1\leq k\leq n}k\left(\left|\hat{f}(n+k)\right|
+\left|\hat{f}(-n-k)\right|\right)\right)$$
holds uniformly for any $x\in [0,\pi]$.

\vspace{3mm}

\bf Proof. \rm The cases $x=0$ and $x=\pi$ are trivial. Suppose $0<x<\pi$. When $0<x\leq\pi/n$, by the inequality $|\sin x|\leq |x|$ we have
$$\left|\sum_{k=1}^{n}\hat{f}(n+k)\sin kx\right|\leq\frac{\pi}{n}\sum_{k=1}^{n}k|\hat{f}(n+k)|\leq\pi\max_{1\leq k\leq n}k\left|\hat{f}(n+k)\right|.\hf (6)$$
On the other hand, in case $\pi/n<x$, we find an natural number $m<n$ such that $m\leq\pi/x<m+1$, then
$$\left|\sum_{k=1}^{n}\hat{f}(n+k)\sin kx\right|\leq\frac{\pi}{m}\sum_{k=1}^{m}k\left|\hat{f}(n+k)\right|+\left|\sum_{k=m+1}^{n}\hat{f}(n+k)\sin kx\right|=:|I_{1}|+|I_{2}|.$$
It is clear that
$$|I_{1}|\leq\pi\max_{1\leq k\leq m}k|\hat{f}(n+k)|
.\hf (7)$$
By Abel transformation,
$$I_{2}=\sum_{k=m+1}^{n-N_{0}-1}\Delta\hat{f}(n+k)\sum_{j=1}^{k}\sin jx+\hat{f}(2n-N_{0})\sum_{j=1}^{n-N_{0}}\sin jx-\hat{f}(n+m+1)\sum_{j=1}^{m}\sin jx$$
$$\hspace{2.0in}+\sum_{k=n-N_{0}+1}^{n}\hat{f}(n+k)\sin kx,$$
so that\footnote{Note that $\sum_{j=1}^{k}\sin jx=O(x^{-1})=O((m+1))$ for $m\leq\pi/x<m+1$.}
$$|I_{2}|\leq (m+1)\sum_{k=m+1}^{n-N_{0}-1}|\Delta\hat{f}(n+k)|+C\max_{m<k\leq n}k|\hat{f}(n+k)|+N_{0}\max_{n-N_{0}<k\leq n}|\hat{f}(n+k)|.$$
But by the condition of Lemma 6, taking a natural number $l$ such that $2^{l}(m+1)\leq n-N_{0}-1<2^{l+1}(m+1)$, and setting
$$\max\limits_{2^{j}(m+1)\leq k<2^{j}(m+1)+N_{0}}|\hat{f}(n+k)|= |\hat{f}(n+k_{j})|,$$
 we then have\footnote{Note that any $k_{j}\leq n$.}
$$
\sum_{k=m+1}^{n-N_{0}-1}|\Delta\hat{f}(n+k)|\leq\sum_{j=0}^{l+1}\sum_{k=2^{j}(m+1)}^{2^{j+1}(m+1)-1}|\Delta\hat{f}(n+k)|\leq M(\mbox{\bf $\hat{f}$})\sum_{j=0}^{l+1}|\hat{f}(n+k_{j})|$$
$$\leq C\max_{1\leq k\leq n}k|\hat{f}(n+k)|\sum_{j=0}^{l+1}k_{j}^{-1}\leq C(m+1)^{-1}\max_{1\leq k\leq n}k|\hat{f}(n+k)|\sum_{j=0}^{\infty}2^{-j}$$
$$\leq C(m+1)^{-1}\max_{1\leq k\leq n}k|\hat{f}(n+k)|.$$
Therefore, with the above estimates, we get
$$|I_{2}|
\leq C\max_{1<k\leq n}k|\hat{f}(n+k)|.$$
Combining with (6) and (7), we have
$$\left|\sum_{k=1}^{n}\hat{f}(n+k)\sin kx\right|=O\left(\max_{1\leq k\leq n}k\left|\hat{f}(n+k)\right|\right).$$

Now write
$$\sum_{k=1}^{n}\hat{f}(-n-k)\sin kx=
\sum_{k=1}^{n}\left(\hat{f}(n+k)+\hat{f}(-n-k)\right)\sin kx-
\sum_{k=1}^{n}\hat{f}(n+k)\sin kx.$$
Applying the above known estimate, by noting that both
$\{\hat{f}(n)\}_{n=0}^{\infty}\in \mbox{\rm GBV}$ and $\{\hat{f}(n)+\hat{f}(-n)\}_{n=1}^{\infty}\in \mbox{\rm GBV}$, we get
$$\left|\sum_{k=1}^{n}\hat{f}(-n-k))\sin kx\right|=
O\left(\max_{1\leq k\leq n}k\left|\hat{f}(n+k)+\hat{f}(-n-k)\right|
+\max_{1\leq k\leq n}k\left|\hat{f}(n+k)\right|\right)$$
$$=O\left(\max_{1\leq k\leq n}k\left(\left|\hat{f}(n+k)\right|+\left|\hat{f}(-n-k)\right|
\right)\right).$$ Lemma 6 is proved.

\vspace{3mm}

\rm Similarly, we can establish the following

\vspace{3mm}

\bf Lemma 7. \it
Let $f\in C_{2\pi}$,
$\{\hat{f}(n)\}_{n=0}^{\infty}\in \mbox{\rm GBV}$ and $\{\hat{f}(n)+\hat{f}(-n)\}_{n=1}^{\infty}\in \mbox{\rm GBV}$,
then
$$\left|\sum_{k=m}^{\infty}\hat{f}(\pm k)\sin kx\right|=
O\left(\max_{k\geq m}k\left(\left|\hat{f}(k)\right|
+\left|\hat{f}(-k)\right|\right)\right)$$
holds uniformly for any $x\in [0,\pi]$ and any $m\geq 1$.

\vspace{3mm}

\bf Proof. \rm Write
$$I(x)=\sum_{k=m}^{\infty}\hat{f}(k)\sin kx.$$
Noting for $x=0$ and $x=\pi$ that
$I(x) = 0$,
we may restrict $x$ within $(0,\pi)$ without loss. Take $N=[1/x]$ and set\footnote{When $N\leq m$, the same argument as in estimating $J_{2}$ can be applied to deal with $I(x)=\sum_{k=m}^{\infty}\hat{f}(k)\sin kx$ directly.}
$$I(x) = \sum_{k=m}^{N-1}\hat{f}(k)\sin kx
+\sum_{k=N}^{\infty}\hat{f}(k)\sin kx =: J_{1}(x) + J_{2}(x).$$
Now
$$|J_{1}(x)| \leq
 \sum_{k=m}^{N-1}|\hat{f}(k)||\sin kx|\leq x\sum_{k=m}^{N-1}k|\hat{f}(k)|$$
 $$< x(N-1)\epsilon_{m} \leq \epsilon_{m}$$
 follows from $N=[1/x]$, where $\epsilon_{m}=\max\limits_{k\geq m}k|\hat{f}(k)|$. By Abel's transformation, similar to the proof of Lemma 6,
$$|J_{2}(x)| = \left|\sum_{k=N}^{\infty}\Delta\hat{f}(k)\sum_{v=1}^{k}\sin vx-\hat{f}(N)\sum_{v=1}^{N-1}\sin vx\right|$$
 $$\leq \sum_{k=N}^{\infty}\left|\Delta\hat{f}(k)\right|\left|\sum_{v=1}^{k}\sin vx\right|+\left|\hat{f}(N)\right|\left|\sum_{v=1}^{N-1}\sin vx\right|$$
$$\leq Cx^{-1}\sum_{j=0}^{\infty}\sum_{k=2^{j}N}^{2^{j+1}N-1}|\Delta\hat{f}(k)|+Cx^{-1}\left|\hat{f}(N)\right|\leq C\epsilon_{m}\sum_{j=0}^{\infty}2^{-j}\leq C\epsilon_{m}.$$
 Altogether, we have
$$|I(x)|\leq C\epsilon_{m},$$
that is the required result for $\sum_{k=m}^{\infty}\hat{f}(k)\sin
kx$. Since $\hat{f}(-k)=\hat{f}(k)+\hat{f}(-k)-\hat{f}(k)$, by the
condition, we can deduce the required result for
$\sum_{k=m}^{\infty}\hat{f}(-k)\sin kx$. Lemma 7 is completed.

\vspace{3mm}

\bf Proof of Theorem 1.
\bf Necessity. \rm
Suppose $f\in C_{2\pi}$. From Lemma 4, (3) clearly holds, while by Lemma 5, we see
$$\sum_{k=2n+1}^{\infty}|\hat{f}(k)+\hat{f}(-k)|
=O(E_{n}(f)),$$
thus (2) holds.

\bf Sufficiency. \rm It can be deduced from Zhou and Le [11, Theorem 1] with a minor revision.

Assume $f\in C_{2\pi}$ now.
By Lemma 4 and Lemma 5, we see
$$\max_{k\geq 1}k|\hat{f}(\pm (n+k))|
+\sum_{k=2n+1}^{\infty}|\hat{f}(k)+\hat{f}(-k)|
=O(E_{n}(f)).$$
On the other hand, rewrite $f(x)$ as
$$f(x)=
\sum_{k=-2n}^{2n}\hat{f}(k)e^{ikx}+i\sum_{k=2n+1}^{\infty}\left(\hat{f}(k)-\hat{f}(-k)\right)\sin kx$$
$$+\frac{1}{2}\sum_{k=2n+1}^{\infty}\left(\hat{f}(k)+\hat{f}(-k)\right)(e^{ikx}+e^{-ikx}),$$
then we have
$$E_{n}(f)\leq \left\|\sum_{k=1}^{n}\left(\hat{f}(n+k)e^{i(n+k)x}
+\hat{f}(-n-k)e^{-i(n+k)x}\right)\right.$$
$$-\left.\sum_{k=1}^{n}\left(\hat{f}(n+k)e^{i(n-k)x}
+\hat{f}(-n-k)e^{i(-n+k)x}\right)\right\|$$
$$+\left\|\sum_{k=2n+1}^{\infty}\left(\hat{f}(k)-\hat{f}(-k)\right)\sin kx\right\|
+\sum_{k=2n+1}^{\infty}\left|\hat{f}(k)+\hat{f}(-k)\right|$$
$$\leq
\left\|\sum_{k=1}^{n}\left(\hat{f}(n+k)\sin kx+\hat{f}(-n-k)\sin(-kx)\right)\right\|
+\left\|\sum_{k=2n+1}^{\infty}\left(\hat{f}(k)-\hat{f}(-k)\right)\sin kx\right\|$$
$$+\sum_{k=2n+1}^{\infty}\left|\hat{f}(k)+\hat{f}(-k)\right|.$$
Applying Lemma 6 immediately yields that
$$\left\|\sum_{k=1}^{n}\left(\hat{f}(n+k)\sin kx+\hat{f}(-n-k)\sin(-kx)\right)\right\|
\leq\max_{1\leq k\leq n}k\left(\left|\hat{f}(n+k)\right|
+\left|\hat{f}(-n-k)\right|\right),$$
while by noting that
$\hat{f}(k)-\hat{f}(-k)=
2\hat{f}(k)-
(\hat{f}(k)+\hat{f}(-k))$ and
both $\{\hat{f}(n)\}_{n=0}^{\infty}\in \mbox{\rm GBV}$ and $\{\hat{f}(n)+\hat{f}(-n)\}_{n=1}^{\infty}\in \mbox{\rm GBV}$, we deduce that
$$\left\|\sum_{k=2n+1}^{\infty}\left(\hat{f}(k)-\hat{f}(-k)\right)\sin kx\right\|=O\left(
2\max_{k\geq 2n+1}k\left|\hat{f}(k) \right| +\max_{k\geq
2n+1}k\left|\hat{f}(k) +\hat{f}(-k)\right|\right)$$ by Lemma 7.
Suppose $\max\limits_{k\geq
2n+1}k\left|\hat{f}(k)\right|=m_{0}\left|\hat{f}(m_{0})\right|$.
Take sufficient large $n$ such that $4n-2N_{0}\geq 3n+1$. Assume
that $2n+1\leq m_{0}\leq 4n$.  Then in view of (1), by a standard
technique used in the present paper, one has
$$\left|\hat{f}(m_{0})\right|\leq \sum_{k=m_{0}}^{4n-2N_{0}}\left|\Delta\hat{f}(k)\right|+\left|\hat{f}(4n-2N_{0}+1)\right|$$
$$\leq \sum_{k=2n-N_{0}}^{4n-2N_{0}}\left|\Delta\hat{f}(k)\right|+\left|\hat{f}(4n-2N_{0}+1)\right|$$
$$\leq C\left(\max_{n-N_{0}\leq k\leq n}\left|\hat{f}(n+k)\right|+\left|\hat{f}(4n-2N_{0}+1)\right|\right).$$
Hence,
$$m_{0}\left|\hat{f}(m_{0})\right|\leq C\left(\max_{1\leq k\leq n}k\left|\hat{f}(n+k)\right|+(4n-2N_{0}+1)\left|\hat{f}(4n-2N_{0}+1)\right|\right).\hf (8)$$
Meanwhile, for $m_{0}\geq 4n+1$,
$$m_{0}\left|\hat{f}(m_{0})
\right|\leq \frac{1}{2}m_{0}\left|\hat{f}(m_{0})-\hat{f}(-m_{0})\right|+\frac{1}{2}m_{0}\left|\hat{f}(m_{0})+\hat{f}(-m_{0})\right|$$
$$\leq
\max_{k\geq 2n+1}k\left|\hat{f}(k)
-\hat{f}(-k)\right|+O\left(\sum_{k=2n+1}^{\infty}
\left|\hat{f}(k)+\hat{f}(-k)\right|\right),$$
where the last inequality follows from $m_{0}\geq 4n+1$ and
the following calculation: Assume that $$\max_{(m_{0}+1)/2+jN_{0}\leq k<(m_{0}+1)/2+(j+1)N_{0}}|\hat{f}(k)+ \hat{f}(-k)|=| \hat{f}(k_{j})+ \hat{f}(-k_{j})|,$$
$$\hspace{3.2in}j=0, 1, \cdots, [(m_{0}+1)/(2N_{0})]-1,$$
$$\hspace{1.0in}2n+1\leq (m_{0}+1)/2+jN_{0}\leq k_{j}< (m_{0}+1)/2+(j+1)N_{0}<m_{0}+1.$$
From condition $(1)$, we get for $0\leq j\leq [(m_{0}+1)/(2N_{0})]-1$,
$$|\hat{f}(m_{0})+\hat{f}(-m_{0})|\leq \sum_{(m_{0}+1)/2+jN_{0}\leq k\leq m_{0}+1+2jN_{0}}|\Delta(\hat{f}(k)+ \hat{f}(-k))|$$
$$\hspace{2.0in}+|\hat{f}(m_{0}+1+2jN_{0}+1)+\hat{f}(-m_{0}-1-2jN_{0}-1)|$$
$$\leq M(\mbox{\bf $\hat{f}$})|\hat{f}(k_{j})+\hat{f}(-k_{j})|+|\hat{f}(m_{0}+1+2jN_{0}+1)+\hat{f}(-m_{0}-1-2jN_{0}-1)|,$$
thus
$$m_{0}|\hat{f}(m_{0})+\hat{f}(-m_{0})|\leq CM(\mbox{\bf $\hat{f}$})|\sum_{k=2n+1}^{\infty}
\left|\hat{f}(k)+\hat{f}(-k)\right|.$$ By the same technique, for
the factor appearing in (8), we also have (note $4n-2N_{0}\geq
3n+1$)
$$(4n-2N_{0}+1)\left|\hat{f}(4n-2N_{0}+1)\right| \leq
\max_{k\geq 2n+1}k\left|\hat{f}(k)
-\hat{f}(-k)\right|+O\left(\sum_{k=2n+1}^{\infty}
\left|\hat{f}(k)+\hat{f}(-k)\right|\right).$$
 From the above discussion, we see that
$$\max_{k\geq 2n+1}k\left|\hat{f}(k)
\right|
\leq C\max_{1\leq k\leq n}k\left|\hat{f}(n+k)\right|
+\max_{k\geq 2n+1}k\left|\hat{f}(k)
-\hat{f}(-k)\right|
+O\left(\sum_{k=2n+1}^{\infty}|\hat{f}(k)+\hat{f}(-k)|\right)$$
holds in any case. With condition that $
\hat{f}(k)+\hat{f}(-k)\in \mbox{\rm GBV}$,
 by a similar way we can also get
$$\max_{k\geq 2n+1}k\left|\hat{f}(k)
+\hat{f}(-k)\right|
\leq C\max_{1\leq k\leq n}k\left(|\hat{f}(n+k)|+|\hat{f}(-n-k)|\right)$$
$$\hspace{2.8in}+O\left(\sum_{k=2n+1}^{\infty}|\hat{f}(k)+\hat{f}(-k)|\right).$$
 Combining all the above estimates, we have completed the proof of Theorem 1.

\vspace{3mm}

\begin{center}
\large\bf \S 3. An Application
\end{center}

\rm
Let
$$f(x) = \sum_{n=0}^{\infty}\hat{f}(n)\cos nx.$$
\rm The following theorem is an interesting application to a hard problem in classical Fourier analysis. Interested readers could check Zhou [10] or Xie and Zhou [9, p.112], for example.

\vspace{3mm}

\bf Theorem 3. \it Let $\{\hat{f}(n)\}\in\mbox{\rm GBV}$ be a real
sequence. If $f\in C_{2\pi}$ and
$$\sum_{k=n+1}^{2n}\hat{f}(k)=O\left(\max_{1\leq k\leq n}k\hat{f}(n+k)\right),$$
then
$$\|f-S_{n}(f)\| = O(E_{n}(f)),$$
where $S_{n}(f)$ is the $n$th Fourier partial sum of $f$.

\begin{center}
{\Large\bf References}
\end{center}

\rm
\begin{enumerate}
\item \rm A. S. Belov, \it On sequential estimate of best approximations and moduli of continuity by sums of trigonometric series with quasimonotone coefficients, \rm Matem. Zametki, 51:4(1992), 132-134, in Russian.
\item L. Leindler, \it On the uniform convergence and boundedness of a certain class of sine series, \rm Anal. Math., 27(2001), 279-285.
\item J. R. Nurcombe, \it On the uniform convergence of sine series with quasimonotone coefficients, \rm J. Math. Anal. Appl., 166(1992), 577-581.
\item V. B. Stanojevic, \it $L^{1}$-convergence of Fourier series with complex quasimonotone coefficients, \rm  Proc. Amer. Math. Soc., 86(1982), 241-247.
\item V. B. Stanojevic, \it $L^{1}$-convergence of Fourier series with $O$-regularly varying quasimonotone coefficients, \rm J. Approx. Theory, 60(1990), 168-173.
\item W. Xiao, T. F. Xie and S. P. Zhou, \it The best approximation rate of certain trigonometric series, \rm Ann. Math. Sinica, 21A(2000), 81-88, in Chinese.
 \item T. F. Xie and S. P. Zhou, \it $L^{1}$-approximation of Fourier series of complex valued functions, \rm Proc. Royal Soc. Edinburgh, 126A(1996), 343-353.
\item T. F. Xie and S. P. Zhou, \it The uniform convergence of certain trigonometric series, \rm J. Math. Anal. Appl., 181(1994), 171-180.
\item T. F. Xie and S. P. Zhou, \it Approximation Theory of Real Functions, \rm Hangzhou Univ. Press, Hangzhou, 1998.
\item S. P. Zhou, \it What strong monotonicity condition on Fourier coefficients can make the ratio $\|f-S_{n}(f)\|/E_{n}(f)$ bounded$?$, \rm Proc. Amer. Math. Soc., 121(1994), 779-785.
\item R. J. Le and S. P. Zhou, \it A new condition for the uniform convergence of certain trigonometric series, \rm Acta Math. Hungar., 108(2005), 161-169.
\end{enumerate}

\begin{flushleft}
\rm S. P. Zhou:\\
Institute of Mathematics\\
Zhejiang Institute of Science and Technology\\
Hangzhou, Zhejiang 310018  China
\end{flushleft}

\begin{flushleft}
\rm R. J. Le:\\
 Department of Mathematics\\
Ningbo University\\
Ningbo, Zhejiang 315211 China\\
 \bf and \rm Institute of Mathematics\\
Zhejiang Institute of Science and Technology\\
Hangzhou, Zhejiang 310018  China
\end{flushleft}

\begin{center}
\bf Keywords \rm uniform convergence,  best approximation, quasimonotone,  bounded variation
\end{center}

\end{document}